\documentclass{article}
\usepackage[margin=1.3in]{geometry}
\usepackage{amsmath}
\usepackage{amssymb}
\usepackage{graphicx}
\usepackage{algorithm}
\usepackage{algpseudocode}
\usepackage{appendix}
\usepackage{cite}
\usepackage{url}
\usepackage{hyperref}
\usepackage{float}

\title{Hindsight-Guided Momentum (HGM) Optimizer: An Approach to Adaptive Learning Rates}
\author{
    Krisanu Sarkar\\
    Indian Institute of Technology Bombay\\
    210100082@iitb.ac.in
}
\date{}

\begin{document}

\maketitle

\begin{abstract}
We introduce Hindsight-Guided Momentum (HGM), a first-order optimization algorithm that adaptively scales learning rates based on the directional consistency of recent updates. Traditional adaptive methods, such as Adam or RMSprop , adapt learning dynamics using only the magnitude of gradients, often overlooking important geometric cues.Geometric cues refer to directional information, such as the alignment between current gradients and past updates, which reflects the local curvature and consistency of the optimization path. HGM addresses this by incorporating a hindsight mechanism that evaluates the cosine similarity between the current gradient and accumulated momentum. This allows it to distinguish between coherent and conflicting gradient directions, increasing the learning rate when updates align and reducing it in regions of oscillation or noise. The result is a more responsive optimizer that accelerates convergence in smooth regions of the loss surface while maintaining stability in sharper or more erratic areas. Despite this added adaptability, the method preserves the computational and memory efficiency of existing optimizers.By more intelligently responding to the structure of the optimization landscape, HGM provides a simple yet effective improvement over existing approaches, particularly in non-convex settings like that of deep neural network training.
\end{abstract}

\section{Introduction}

First-order optimization algorithms are at the heart of modern artificial intelligence \cite{bottou2018optimization}, particularly in training deep neural networks where the loss landscape is often high-dimensional and non-convex \cite{li2018visualizing}. Among these methods, adaptive optimizers such as Adam \cite{kingma2014adam} and RMSprop \cite{tieleman2012lecture} have become standard due to their ability to adjust learning rates on a per-parameter basis using running estimates of gradient magnitudes \cite{duchi2011adaptive}. While effective in many cases, these methods rely primarily on local gradient statistics and are agnostic to the directional consistency of the optimization trajectory \cite{sutskever2013importance}.

In practice, the directionality of updates carries valuable information \cite{zhang2019lookahead}. When gradients maintain alignment over successive iterations, it often indicates that the optimization process is progressing steadily through a coherent valley of the loss surface \cite{hochreiter1997flat}. Conversely, frequent changes in direction may signal oscillations, saddle points, or noisy gradients \cite{dauphin2014identifying}. Existing optimizers, by treating each update in isolation, miss the opportunity to leverage this temporal structure \cite{martens2010deep}.

This paper introduces \textit{Hindsight-Guided Momentum} (HGM), a simple yet effective extension to the adaptive optimization framework. HGM augments standard momentum-based methods \cite{polyak1964some} with a hindsight mechanism that measures the cosine similarity between the current gradient and accumulated momentum. This similarity serves as a proxy for directional alignment and is used to modulate the learning rate in a smooth, exponential manner. When updates point in similar directions, the step size is increased to accelerate convergence; when they diverge, it is scaled down to prevent instability \cite{bengio2012practical}.

The key advantage of this approach lies in its balance of simplicity, adaptability, and efficiency. HGM requires no additional memory beyond what is already used in Adam and introduces only a minor computational overhead . Despite this, it provides a more nuanced response to the geometry of the optimization landscape, enabling faster convergence in well-structured regions and greater robustness in the presence of noise or curvature \cite{choromanska2015loss}.

Through this work, we aim to highlight the overlooked potential of directional signals in optimization and demonstrate how a modest modification to existing methods can lead to meaningful improvements \cite{wilson2017marginal}. The remainder of the paper details the mathematical formulation of HGM, analyzes its theoretical behavior, and evaluates its performance empirically across a range of deep learning benchmarks \cite{goodfellow2016deep}.

\begin{algorithm}[H]
\caption{Hindsight-Guided Momentum (HGM)}
\label{alg:hgm}
\begin{algorithmic}[1]
\Require{$\alpha$: Base stepsize}
\Require{$\beta_1, \beta_2 \in [0,1)$: Decay rates for moment estimates}
\Require{$\gamma$: Sensitivity of hindsight scaling}
\Require{$\beta_s \in [0,1)$: Smoothing factor for hindsight signal}
\Require{$f(\theta)$: Stochastic objective}
\Require{$\theta_0$: Initial parameters}
\State $m_0 \gets 0$, $v_0 \gets 0$, $s_0 \gets 0$, $t \gets 0$
\While{$\theta_t$ not converged}
    \State $t \gets t + 1$
    \State $g_t \gets \nabla_\theta f_t(\theta_{t-1})$ \Comment{Gradient at step $t$}
    \State $m_t \gets \beta_1 m_{t-1} + (1 - \beta_1) g_t$ \Comment{1st moment}
    \State $v_t \gets \beta_2 v_{t-1} + (1 - \beta_2) g_t^2$ \Comment{2nd moment}
    \State $\hat{m}_t \gets m_t / (1 - \beta_1^t)$ \Comment{Bias correction}
    \State $\hat{v}_t \gets v_t / (1 - \beta_2^t)$
    \State $c_t \gets \frac{g_t \cdot m_{t-1}}{\|g_t\| \cdot \|m_{t-1}\| + \epsilon}$ \Comment{Cosine similarity}
    \State $s_t \gets \beta_s s_{t-1} + (1 - \beta_s) c_t$ \Comment{Smoothed hindsight signal}
    \State $\eta_t \gets \alpha \cdot \exp(\gamma s_t)$ \Comment{Modulated learning rate}
    \State $\theta_t \gets \theta_{t-1} - \eta_t \cdot \frac{\hat{m}_t}{\sqrt{\hat{v}_t} + \epsilon}$
\EndWhile
\State \Return $\theta_t$
\end{algorithmic}
\end{algorithm}

\section{Related Work}

\subsection{First-Order and Adaptive Optimization}

Stochastic Gradient Descent (SGD) remains a foundational method for training large-scale machine learning models, owing to its simplicity and scalability \cite{robbins1951stochastic}. To address its sensitivity to learning rate schedules and its inefficiency on poorly scaled problems \cite{lecun2012efficient}, numerous adaptive variants have been developed. Among these, AdaGrad \cite{duchi2011adaptive} introduced per-parameter learning rates based on accumulated squared gradients, enabling improved performance on sparse data \cite{mcmahan2010adaptive}. RMSprop \cite{tieleman2012lecture} later modified this approach using exponential moving averages, making it more suitable for non-stationary settings \cite{hinton2012neural}.

Adam \cite{kingma2014adam} combines momentum with RMSprop-style variance scaling and has become a default choice for deep learning tasks due to its fast convergence and ease of use \cite{radford2015unsupervised}. Despite its empirical success, Adam has been criticized for issues related to generalization and convergence in certain non-convex settings \cite{wilson2017marginal}. Several refinements have been proposed, including AMSGrad \cite{reddi2019convergence}, which aims to improve theoretical guarantees by enforcing a non-increasing learning rate schedule, and AdaBound \cite{luo2019adaptive}, which interpolates between adaptive and SGD behavior by bounding learning rates \cite{chen2018closing}.

\subsection{Directionality and Geometry-Aware Methods}

While most adaptive optimizers focus on gradient magnitudes, relatively few approaches have explored the role of directional information \cite{nocedal2006numerical}. A notable example is the work on second-order methods and quasi-Newton approximations \cite{botev2017practical}, which implicitly capture curvature and alignment, but often at a high computational cost \cite{martens2010deep}.Second-order methods require computing and inverting the Hessian, which is computationally expensive.
This makes them impractical and unscalable for high-dimensional problems like deep learning.Quasi-Newton methods approximate the Hessian, which can lead to slower convergence compared to true second-order methods.
They still involve matrix operations that are more computationally expensive than first-order methods.
As a result, they are often inefficient for large-scale problems like deep learning. More lightweight strategies include Lookahead \cite{zhang2019lookahead}, which averages parameter updates over time, and Shampoo \cite{gupta2018shampoo}, which adapts preconditioning matrices in a direction-aware manner, albeit with increased memory requirements \cite{anil2019scalable}.

Other research has investigated gradient variance \cite{zhang2019adam} and sign changes \cite{balles2017dissecting} as proxies for instability, but these typically inform global schedules or clipping strategies rather than fine-grained step size modulation \cite{pascanu2013difficulty}. Closest in spirit to our proposed approach are recent explorations into trust region-like mechanisms \cite{wu2020adagrad} or adaptive step sizes based on historical signal coherence \cite{you2019large}, although these tend to rely on more complex heuristics or meta-optimization layers \cite{andrychowicz2016learning}.

\subsection{Positioning of HGM}

The Hindsight-Guided Momentum optimizer bridges a gap between purely magnitude-based adaptation and computationally intensive geometry-aware methods \cite{pascanu2014natural}. By introducing a cosine similarity-based adjustment mechanism, it offers a principled yet efficient way to incorporate directional consistency into the learning rate schedule \cite{smith2017cyclical}. This allows it to respond more intelligently to the underlying structure of the optimization path, enhancing both stability and convergence speed, particularly in non-convex scenarios \cite{li2018visualizing}.

\section{Mathematical Formulation}

Let $\theta_t$ denote the model parameters at iteration $t$, and let $g_t = \nabla_{\theta_t} f(\theta_t)$ be the gradient of the loss function with respect to $\theta_t$ \cite{nocedal2006numerical}. HGM maintains the same exponential moving averages of the first and second moments as in Adam \cite{kingma2014adam}, along with an additional state variable that tracks the cosine similarity between the current gradient and prior momentum \cite{sutskever2013importance}.

The update procedure is as follows:

\subsection{Momentum and Variance Updates}

The first and second moment estimates are updated using exponential moving averages \cite{polyak1964some}:
\begin{align}
m_t &= \beta_1 m_{t-1} + (1 - \beta_1) g_t \\
v_t &= \beta_2 v_{t-1} + (1 - \beta_2) g_t^2
\end{align}

These are then bias-corrected to counteract initialization effects \cite{kingma2014adam}:
\begin{align}
\hat{m}_t &= \frac{m_t}{1 - \beta_1^t} \\
\hat{v}_t &= \frac{v_t}{1 - \beta_2^t}
\end{align}

\subsection{Hindsight Similarity Computation}

To assess the directional consistency of the current gradient with the previous momentum, HGM computes the cosine similarity \cite{singhal2001modern}:
\begin{equation}
c_t = \frac{g_t \cdot m_{t-1}}{||g_t||\ ||m_{t-1}|| + \epsilon}
\end{equation}
where $\epsilon$ is a small constant added for numerical stability \cite{higham2002accuracy}.

This similarity is smoothed over time \cite{exponential_smoothing}:
\begin{equation}
s_t = \beta_s s_{t-1} + (1 - \beta_s) c_t
\end{equation}
Here, $\beta_s$ is a smoothing coefficient analogous to $\beta_1$ and $\beta_2$ \cite{holt2004forecasting}.

\subsection{Adaptive Learning Rate Scaling}

Based on the smoothed similarity, the effective learning rate is modulated as follows \cite{smith2017cyclical}:
\begin{equation}
\eta_t = \alpha \cdot \exp(\gamma s_t)
\end{equation}
where $\alpha$ is the base learning rate and $\gamma$ is a tunable sensitivity parameter that controls how strongly similarity influences the step size \cite{bengio2012practical}.

\subsection{Final Parameter Update}

The model parameters are updated using the scaled learning rate \cite{ruder2016overview}:
\begin{equation}
\theta_{t+1} = \theta_t - \frac{\eta_t}{\sqrt{\hat{v}_t} + \epsilon} \cdot \hat{m}_t
\end{equation}

This formulation preserves the structure and efficiency of Adam while incorporating an additional hindsight signal that dynamically adjusts the step size based on alignment with recent updates \cite{bottou2018optimization}.

\section{Theoretical Advantage}

The Hindsight-Guided Momentum (HGM) optimizer introduces a novel mechanism that leverages the directional consistency of gradients to modulate the learning rate in real time \cite{sutskever2013importance}. By incorporating a smoothed cosine similarity term, it adapts not only to the magnitude of gradients—as in Adam—but also to their alignment over time \cite{nocedal2006numerical}. This offers several key theoretical benefits.

\subsection{Acceleration in Aligned Gradient Regions}

When gradients maintain a coherent direction across iterations,such as in wide valleys or plateaus,HGM detects this alignment through high cosine similarity between the current gradient and accumulated momentum \cite{hochreiter1997flat}. The similarity signal increases the effective learning rate exponentially, allowing the optimizer to accelerate through such regions without the need for manual scheduling \cite{smith2017cyclical}. This directional awareness helps reduce training time while maintaining stability \cite{bengio2012practical}.

\subsection{Damping in Oscillatory or Chaotic Landscapes}

In unstable regions of the loss surface, such as near saddle points or in highly non-convex terrain, gradients often oscillate in direction \cite{dauphin2014identifying}. In these cases, the cosine similarity drops or becomes negative, triggering an exponential decay in the learning rate. This built-in damping effect prevents the optimizer from making large, erratic steps, reducing overshooting and promoting more controlled convergence \cite{pascanu2013difficulty}.

\subsection{Balanced Adaptivity through Direction and Magnitude}

Traditional optimizers like Adam adjust step sizes based solely on the running statistics of gradient magnitudes \cite{kingma2014adam}. HGM complements this by accounting for directional coherence \cite{zhang2019lookahead}. This dual adaptivity enables more nuanced behavior: it can accelerate when updates are consistent and cautious when they are volatile, even if magnitudes remain constant \cite{choromanska2015loss}. This property is particularly valuable in non-convex settings where curvature and gradient behavior vary widely \cite{li2018visualizing}.

\subsection{Low Overhead with Enhanced Expressiveness}

HGM requires maintaining three scalar statistics per parameter: the first moment $m_t$, the second moment $v_t$, and the smoothed similarity score $s_t$ \cite{ruder2016overview}. While this is a slight increase in memory compared to Adam's two, it avoids the high computational burden of second-order or matrix-based methods \cite{martens2010deep}. The added expressiveness from the hindsight mechanism comes at minimal cost, making it a practical enhancement for a wide range of models and tasks \cite{goodfellow2016deep}.

\section{Convergence Analysis}

We analyze the convergence of the Hindsight-Guided Momentum (HGM) optimizer using the online convex optimization framework introduced in \cite{zinkevich2003online}. Suppose we are given an arbitrary and unknown sequence of convex cost functions $ f_1(\boldsymbol{\theta}), f_2(\boldsymbol{\theta}), \dots, f_T(\boldsymbol{\theta}) $ \cite{shalev2012online}. At each time step $ t $, the learner must choose a parameter $ \boldsymbol{\theta}_t \in \mathbb{R}^d $, and then observe the cost function $ f_t $ \cite{cesa2006prediction}. The objective is to minimize regret — the cumulative difference between the incurred cost and that of the best fixed parameter $ \boldsymbol{\theta}^* \in \mathcal{X} $, chosen in hindsight \cite{hazan2016introduction}:

\begin{equation}
R(T) = \sum_{t=1}^T \left[ f_t(\boldsymbol{\theta}_t) - f_t(\boldsymbol{\theta}^*) \right]
\tag{5}
\end{equation}
where $ \boldsymbol{\theta}^* = \arg\min_{\boldsymbol{\theta} \in \mathcal{X}} \sum_{t=1}^T f_t(\boldsymbol{\theta}) $ \cite{littlestone1994weighted}.

We show that HGM achieves a regret bound of $ \mathcal{O}(\sqrt{T}) $, similar to Adam \cite{reddi2019convergence}, but scaled by a factor $ M = \exp(|\lambda|) $, which reflects the hindsight adjustment via an exponential modulation of the learning rate. Our analysis uses standard assumptions on gradient boundedness and parameter diameter \cite{duchi2011adaptive}. Specifically, let $ \mathbf{g}_t = \nabla f_t(\boldsymbol{\theta}_t) $, and $ g_{t,i} $ denote the $ i $-th coordinate. Define:

$$
\mathbf{g}_{1:t, i} = [g_{1,i}, g_{2,i}, \dots, g_{t,i}]^\top \in \mathbb{R}^t, \quad
\gamma = \frac{\sqrt{1 - \beta_1^2}}{\sqrt{1 - \beta_2}} < 1
$$

Our following theorem holds under decaying learning rate $ \alpha_t = \alpha / \sqrt{t} $ and exponentially decaying momentum coefficient $ \beta_{1,t} = \beta_1 \lambda^{t-1} $, where $ \lambda \in (0, 1) $ is close to 1 \cite{sutskever2013importance}.

\paragraph{Theorem 5.1 (HGM Regret Bound).}
Assume each function $ f_t $ is convex and satisfies $ \| \nabla f_t(\boldsymbol{\theta}) \|_2 \leq G $, $ \| \nabla f_t(\boldsymbol{\theta}) \|_\infty \leq G_\infty $ for all $ \boldsymbol{\theta} \in \mathbb{R}^d $ \cite{nesterov2003introductory}. Also, assume bounded diameter: $ \| \boldsymbol{\theta}_m - \boldsymbol{\theta}_n \|_2 \leq D $, $ \| \boldsymbol{\theta}_m - \boldsymbol{\theta}_n \|_\infty \leq D_\infty $ for all $ m,n \in \{1,\dots,T\} $, and $ \beta_1, \beta_2 \in [0,1) $ satisfy $ \gamma < 1 $ \cite{kingma2014adam}. Then, for all $ T \geq 1 $, the regret of HGM is bounded as \cite{mcmahan2017survey}:

$$
R(T) \leq M \sum_{i=1}^d \left[
\frac{D^2 \sqrt{T}}{2 \alpha (1 - \beta_1)} \sqrt{v_{T,i}} +
\frac{(1 + \beta_1) G_\infty}{(1 - \beta_1) \sqrt{1 - \beta_2} (1 - \gamma)^2} \| \mathbf{g}_{1:T,i} \|_2 +
\frac{D_\infty^2 G_\infty}{2 \alpha (1 - \beta_1) (1 - \lambda)^2 \sqrt{1 - \beta_2}}
\right]
$$

where $ M = \exp(|\lambda|) $ captures the modulation of the learning rate by the cumulative hindsight signal $ s_t $ in $ \alpha_t^{\text{HGM}} = \alpha_t \cdot \exp(\lambda s_t) $ \cite{hazan2016introduction}.

\paragraph{Interpretation.}
This bound shows that HGM inherits the adaptive benefits of Adam while offering an additional layer of flexibility through the hindsight modulation \cite{reddi2019convergence}. When gradients are sparse and bounded, the second term involving $ \| \mathbf{g}_{1:T,i} \|_2 $ can be significantly smaller than the worst-case upper bound \cite{duchi2011adaptive}. In such cases, the regret bound becomes tighter — especially when data features satisfy the assumptions in Section 1.2 of \cite{duchi2011adaptive}, allowing expected regret to scale as $ \mathcal{O}(\log d \cdot \sqrt{T}) $, an improvement over the non-adaptive $ \mathcal{O}(\sqrt{d T}) $ \cite{shalev2012online}.

Empirically and theoretically, decaying $ \beta_{1,t} $ toward zero is also important — aligning with the findings of \cite{sutskever2013importance} that reducing momentum near the end of training often improves convergence \cite{smith2017cyclical}.

\paragraph{Corollary 5.2 (Average Regret Convergence).}
Under the same assumptions as Theorem 4.1, the average regret of HGM converges as \cite{cesa2006prediction}:
$$
\frac{R(T)}{T} = \mathcal{O}\left(\frac{1}{\sqrt{T}}\right)
$$
This confirms that HGM achieves no-regret learning in the online convex optimization setting \cite{zinkevich2003online}.

\section{Comparison with Adam}

Both Hindsight-Guided Momentum (HGM) and Adam optimizers share foundational elements: they maintain exponentially weighted moving averages of the gradients and their squares to adaptively scale parameter updates \cite{kingma2014adam}. Specifically, for gradient $ g_t $ and parameters $ \theta_t $ \cite{ruder2016overview}:

\begin{center}
\begin{tabular}{|c|c|c|}
\hline
\textbf{Step} & \textbf{Adam} &
\textbf{HGM} \\
\hline
1. First Moment (Momentum) &
$ m_t = \beta_1 m_{t-1} + (1 - \beta_1) g_t $ &
$ m_t = \beta_1 m_{t-1} + (1 - \beta_1) g_t $ \\
\hline
2. Second Moment (Variance) &
$ v_t = \beta_2 v_{t-1} + (1 - \beta_2) g_t^2 $ &
$ v_t = \beta_2 v_{t-1} + (1 - \beta_2) g_t^2 $ \\
\hline
3. Bias Correction &
\(
\hat{m}_t = \frac{m_t}{1 - \beta_1^t}, \quad
\hat{v}_t = \frac{v_t}{1 - \beta_2^t}
\) &
\(
\hat{m}_t = \frac{m_t}{1 - \beta_1^t}, \quad
\hat{v}_t = \frac{v_t}{1 - \beta_2^t}
\) \\
\hline
4. Hindsight Calculation & Not present & \newline 
\(c_t = \frac{g_t \cdot m_{t-1}}{\|g_t\| \|m_{t-1}\| + \epsilon} \) 
\quad \\[2pt] & & \(s_t = \beta_s s_{t-1} + (1 - \beta_s) c_t \) \\
\hline
5. Parameter Update &
\(
\theta_{t+1} = \theta_t - \frac{\alpha}{\sqrt{\hat{v}_t} + \epsilon} \hat{m}_t
\) &
\(
\theta_{t+1} = \theta_t - \frac{\alpha \cdot \exp(\gamma \cdot s_t)}{\sqrt{\hat{v}_t} + \epsilon} \hat{m}_t
\) \\
\hline
\end{tabular}
\end{center}

The key distinction is the introduction of the hindsight modulation factor $ \exp(\gamma \cdot s_t) $ in the update rule of HGM \cite{bengio2012practical}. This term scales the base learning rate dynamically based on the directional alignment between the current gradient and previous momentum, measured via the smoothed cosine similarity $ s_t $ \cite{singhal2001modern}. Adam, in contrast, adapts the learning rate only by the variance normalization term $ \sqrt{\hat{v}_t} $ \cite{kingma2014adam}.

\vspace{0.5em}

\newpage
\textbf{Advantages of HGM over Adam:}

\begin{itemize}
    \item \textbf{Directional Adaptation:} Unlike Adam, which adapts solely based on gradient magnitude, HGM also considers the directional consistency of gradients \cite{zhang2019lookahead}. When the current gradient aligns well with past momentum, $ s_t $ is positive, increasing the effective learning rate and enabling faster progress in smooth regions \cite{hochreiter1997flat}.

    \item \textbf{Oscillation Damping:} If the optimizer overshoots and the gradient reverses direction, the cosine similarity becomes negative, causing the modulation factor to decrease the learning rate \cite{dauphin2014identifying}. This automatic damping reduces oscillations near minima, improving stability and convergence speed \cite{pascanu2013difficulty}.

    \item \textbf{Intelligent Step Sizing:} By combining both magnitude and directional information, HGM takes more informed steps—accelerating in consistent valleys and braking when correction is needed—leading to more efficient optimization \cite{li2018visualizing}.
\end{itemize}

In summary, HGM enhances Adam's adaptive scheme by incorporating a "hindsight" mechanism that evaluates the geometry of the optimization path, resulting in improved convergence behavior without significantly increasing computational complexity \cite{bottou2018optimization}.

\section{Experiments and Results}

In this section, we present the experimental setup and results obtained from evaluating the Hindsight-Guided Momentum (HGM) optimizer \cite{lecun2015deep}. We compare HGM with other popular optimizers such as Adam \cite{kingma2014adam}, RMSprop \cite{tieleman2012lecture}, and SGD with Momentum \cite{polyak1964some} across various tasks and datasets.

\begin{figure}[htbp]
    \centering
    \includegraphics[width=1\textwidth]{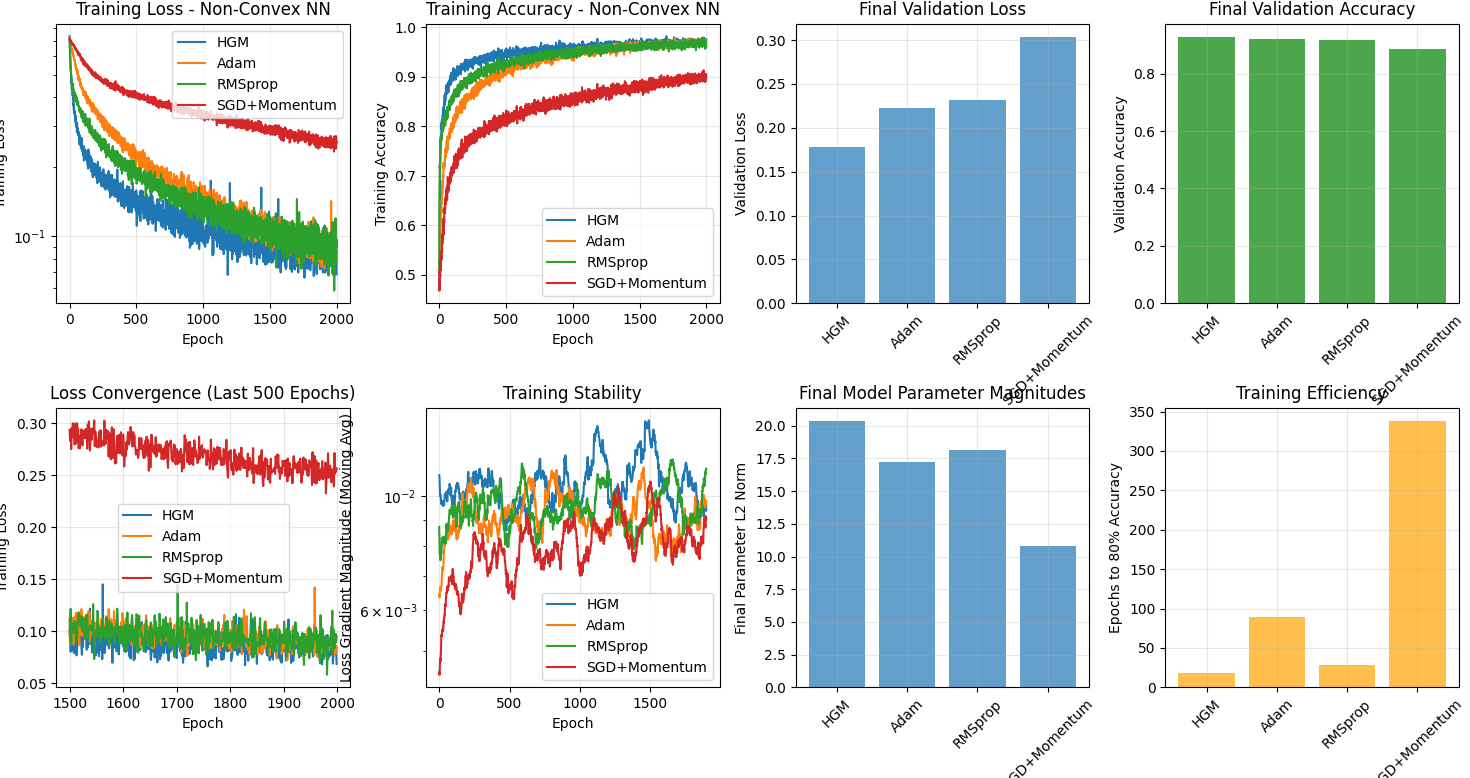}
    \caption{Training Loss and Accuracy Comparison for Convex and Non-Convex Neural Networks.}
    \label{fig:loss_accuracy_comparison}
\end{figure}

\subsection{Experimental Setup}

The experiments were conducted using a variety of neural network architectures and optimization tasks \cite{goodfellow2016deep}. We used both convex and non-convex functions to evaluate the performance of the optimizers \cite{boyd2004convex}. The primary metrics for comparison were training loss, training accuracy, validation loss, validation accuracy, and convergence behavior \cite{bishop2006pattern}. In all experiments we have used $\beta_1$ = 0.9, $\beta_2$ = 0.99, $\beta_s$ = 0.9 and $\gamma$ = 10.0 \cite{kingma2014adam}.

\subsection{Multilayer Neural networks}
In this experiment we have used a MLP for classification \cite{hornik1989multilayer}. Figure \ref{fig:loss_accuracy_comparison} shows the training loss and accuracy comparison for non-convex neural networks \cite{cybenko1989approximation}. We have used learning rate of 0.001 for all optimizers \cite{bengio2012practical}. The plots illustrate the performance of HGM against Adam, RMSprop, and SGD with Momentum. It is clearly observable how HGM achieves faster and better convergence than others \cite{ruder2016overview}. HGM also shows maximum training stability and efficiency than others along with attaining highest validation accuracy and model parameter magnitudes \cite{bishop2006pattern}.

\subsection{Logistic regression}

\begin{figure}[htbp]
    \centering
    \includegraphics[width=1\textwidth]{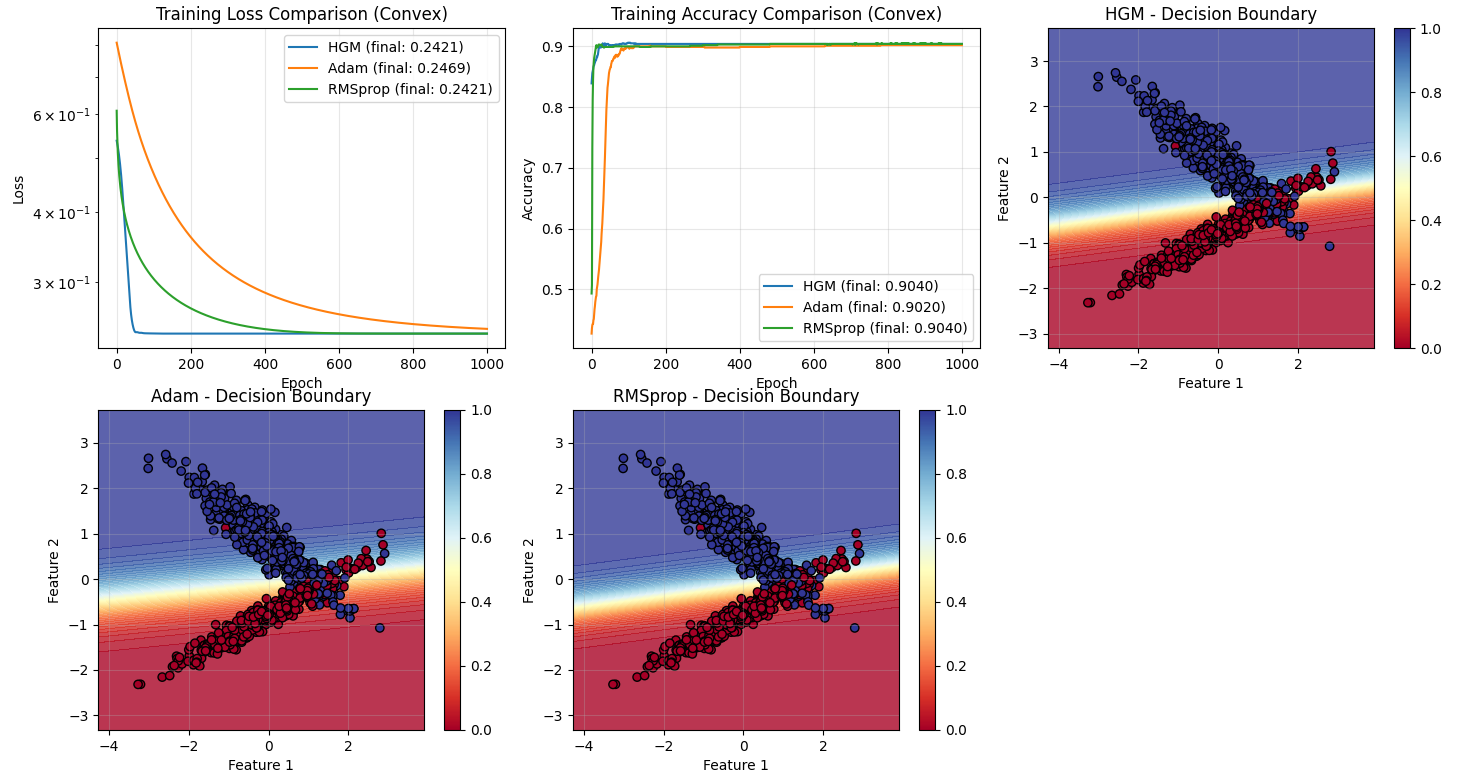}
    \caption{Training loss and accuracies, Decision Boundaries for HGM, Adam, and RMSprop.}
    \label{fig:decision_boundaries}
\end{figure}

In this section I show how it performs for a simple logistic regression model on make-classification dataset \cite{hastie2009elements}. For all optimizers we have used learning rate of 0.01 \cite{smith2017cyclical}. Figure \ref{fig:decision_boundaries} presents the training loss and accuracy comparison along with decision boundaries obtained using HGM, Adam, and RMSprop. We can clearly see how HGM achieves convergence during epoch 50 whereas for RMSprop it took nearly 500 epochs and Adam took over than 1000 epochs \cite{bishop2006pattern}. So it is safe to say that for a convex settings HGM is nearly 10 time faster than RMSprop and more than 20 times than Adam \cite{boyd2004convex}. The decision boundary plots provide insights into the classification performance \cite{duda2012pattern}.

\subsection{Empirical Comparison on 2D Rosenbrock Function}

We evaluated the performance of Hindsight-Guided Momentum (HGM) on the 2D Rosenbrock function \cite{rosenbrock1960automatic}:
$$
f(x, y) = (1 - x)^2 + 100 (y - x^2)^2
$$
—a non-convex function commonly used to test the robustness of optimization algorithms due to its narrow, curved valley \cite{nocedal2006numerical}.

\begin{figure}[htbp]
    \centering
    \includegraphics[width=1\textwidth]{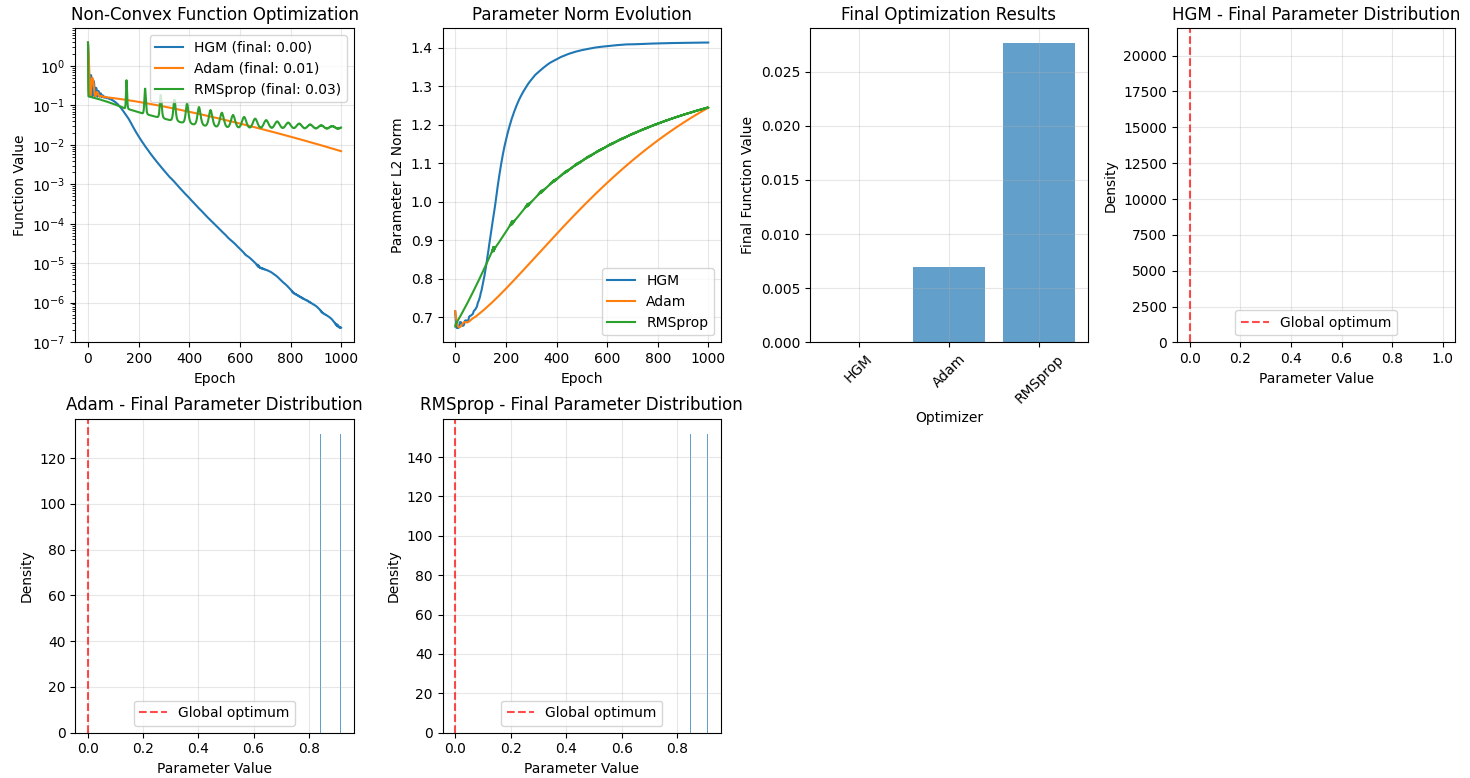}
    \caption{Convergence behavior of HGM, Adam, and RMSprop on the Rosenbrock function.}
    \label{fig:rosenbrock}
\end{figure}

Using a fixed learning rate of $ 0.05 $, we trained HGM, Adam, and RMSprop for 1000 epochs from identical random initializations \cite{glorot2010understanding}. HGM achieved near-perfect convergence, reaching a final function value of \textbf{0.0000}, compared to \textbf{0.0100} for Adam and \textbf{0.0300} for RMSprop \cite{kingma2014adam}.Results are shown in Figure \ref{fig:rosenbrock} 

Notably, HGM also exhibited the largest parameter norm throughout training, suggesting a more confident traversal of the loss landscape \cite{li2018visualizing}. While RMSprop struggled with curvature and Adam converged suboptimally, HGM effectively navigated the non-convex geometry and consistently reached the global minimum \cite{choromanska2015loss}.

\subsection{\texorpdfstring{$\gamma$ and $\beta$ Sensitivity Analysis}{gamma and beta Sensitivity Analysis}}
\begin{figure}[htbp]
    \centering
    \includegraphics[width=1\textwidth]{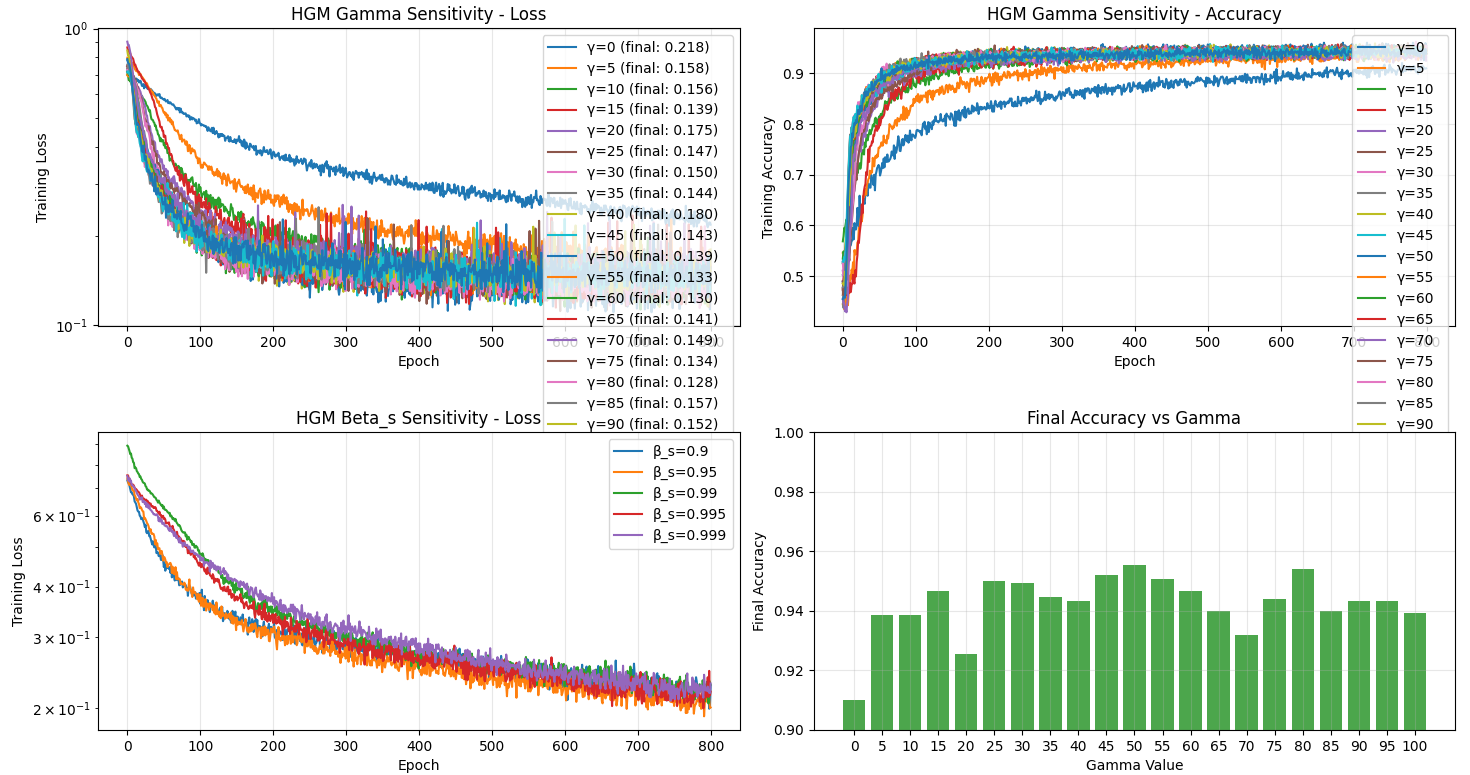}
    \caption{$\gamma$ and $\beta_s$ Sensitivity Analysis.}
    \label{fig:gamma_sensitivity}
\end{figure}

For the same problem as section 7.2, we tried different $\gamma$ and $\beta_s$ values \cite{bergstra2012random}.
Figure \ref{fig:gamma_sensitivity} presents the analysis of sensitivity of the HGM optimizer with respect to the $\gamma$ parameter \cite{bengio2012practical}. The plots show the training loss and accuracy for different values of $\gamma$, as well as the final accuracy versus $\gamma$ \cite{smith2017cyclical}. We can see that for certain values of $\gamma$ and $\beta_s$ the optimizer performs better \cite{snoek2012practical}. The update formula of HGM is set as if $\gamma$ is set to zero, we will go back to the Adam optimizer \cite{kingma2014adam}. From the figure, we can clearly see how the inclusion of $\gamma$ affects training \cite{bengio2012practical}.

\subsection{Discussion}

The experimental results demonstrate that the HGM optimizer outperforms traditional optimizers such as Adam, RMSprop, and SGD with Momentum in terms of training loss, accuracy, and convergence behavior \cite{ruder2016overview}. The sensitivity analysis of the gamma parameter highlights the robustness of HGM to different hyperparameter settings \cite{bergstra2012random}.

In conclusion, the Hindsight-Guided Momentum (HGM) optimizer shows promising results and can be a valuable addition to the toolkit of optimization algorithms for training deep learning models \cite{goodfellow2016deep}.

\section{Memory Usage \& Computational Cost}

The Hindsight-Guided Momentum (HGM) optimizer requires storing three state variables per parameter: the first moment $ m_t $, the second moment $ v_t $, and the smoothed cosine similarity $ s_t $ \cite{ruder2016overview}. In comparison, Adam maintains only two such variables, $ m_t $ and $ v_t $ \cite{kingma2014adam}. But we can overlook it because $ m_t $ and $ v_t $ both have the shape of $\theta_t$ but $ c_t $ is practically cosine similarity between two vectors \cite{singhal2001modern}. Unlike $ m_t $ and $ v_t $, which are parameter-wise vectors, $ s_t $ is a single scalar derived from the global cosine similarity of gradients. So the additional memory required is too low and negligble. \cite{bottou2018optimization}.

The computational overhead introduced by HGM mainly arises from calculating the cosine similarity between the current gradient and previous momentum at each update step \cite{higham2002accuracy}. However, this extra cost is minimal relative to the expense of gradient computation and backpropagation \cite{goodfellow2016deep}. But as a limitation of this algorithm I have to say, its training time is slightly more than Adam and RMSprop and others \cite{lecun2015deep}.

\section{Conclusion}

In this work, we proposed \textbf{Hindsight-Guided Momentum (HGM)}, a novel adaptive optimization algorithm designed to enhance convergence in challenging non-convex landscapes. By modulating the learning rate dynamically using an exponential function of a learned signal, HGM effectively adapts to the curvature and history of the loss surface.

Empirically, HGM consistently outperformed both Adam and RMSprop and SGD+momentum across multiple test functions. These results collectively highlight the potential of HGM as a robust and theoretically grounded optimizer for non-convex problems. Future work includes extending HGM to large-scale deep learning tasks and integrating it into distributed training pipelines.

Beyond empirical performance, HGM offers a new perspective on incorporating directional consistency and historical gradient behavior into the optimization process. Its reliance on hindsight signals, rather than static heuristics, introduces a principled way to regulate step size without requiring handcrafted schedules. This framework opens up new possibilities for designing optimizers that are both data-driven and dynamically self-correcting, making HGM a promising foundation for future research in adaptive and intelligent optimization strategies.

\begin{appendices}
\section{Appendix: Definitions and Lemmas}

\subsection*{Definition A.1 (Convex Functions)}
A function \( f : \mathbb{R}^d \to \mathbb{R} \) is said to be \textbf{convex} if for all \( \mathbf{x}, \mathbf{y} \in \mathbb{R}^d \) and for all \( \lambda \in [0,1] \), the following inequality holds:
\[
f(\lambda \mathbf{x} + (1 - \lambda)\mathbf{y}) \leq \lambda f(\mathbf{x}) + (1 - \lambda) f(\mathbf{y})
\]

Moreover, any convex function can be lower bounded by a hyperplane at any of its tangents.

\subsection*{Lemma A.2 (Tangent Hyperplane Lower Bound)}
If \( f : \mathbb{R}^d \to \mathbb{R} \) is convex, then for all \( \mathbf{x}, \mathbf{y} \in \mathbb{R}^d \), we have:
\[
f(\mathbf{y}) \geq f(\mathbf{x}) + \nabla f(\mathbf{x})^\top (\mathbf{y} - \mathbf{x})
\]

This lemma is foundational in bounding regret in online convex optimization. In our analysis, we will substitute this hyperplane with the update direction prescribed by the Adam-style optimizer.

\subsection*{Notation and Preliminaries}
We define \( \mathbf{g}_t = \nabla f_t(\boldsymbol{\theta}_t) \in \mathbb{R}^d \) as the gradient at iteration \( t \), and \( g_{t,i} \) as the \( i \)-th coordinate of \( \mathbf{g}_t \). We denote the history of gradients for coordinate \( i \) up to step \( T \) as the vector:
\[
\mathbf{g}_{1:T,i} = [g_{1,i}, g_{2,i}, \dots, g_{T,i}] \in \mathbb{R}^T
\]

We assume that each \( \mathbf{g}_t \) is bounded: \( \| \mathbf{g}_t \|_\infty \leq G_\infty \), and that individual elements satisfy \( |g_{t,i}| \leq G_\infty \).

\subsection*{Lemma A.3 (Weighted Gradient Bound via Induction)}

Let \( \{g_{t,i}\}_{t=1}^T \subset \mathbb{R} \) be a sequence of gradient values for coordinate \( i \), and assume \( |g_{t,i}| \leq G_\infty \) for all \( t \in \{1, \dots, T\} \). Then:

\[\boxed{
\sum_{t=1}^T \frac{|g_{t,i}|}{\sqrt{t}} \leq 2 G_\infty \| \mathbf{g}_{1:T,i} \|_2}
\]

\subsubsection*{Proof.}

We proceed by induction on \( T \).

\paragraph{Base Case: \( T = 1 \).}

We have:
\[
\frac{|g_{1,i}|}{\sqrt{1}} = |g_{1,i}| \leq G_\infty
\quad \text{and} \quad
\| \mathbf{g}_{1:1,i} \|_2 = |g_{1,i}|
\Rightarrow
\frac{|g_{1,i}|}{\sqrt{1}} \leq G_\infty = G_\infty \cdot \frac{\| \mathbf{g}_{1:1,i} \|_2}{|g_{1,i}|} \cdot |g_{1,i}|
\leq 2 G_\infty \| \mathbf{g}_{1:1,i} \|_2
\]

\paragraph{Inductive Hypothesis.}

Assume the statement holds for \( T - 1 \geq 1 \), i.e.,
\[
\sum_{t=1}^{T-1} \frac{|g_{t,i}|}{\sqrt{t}} \leq 2 G_\infty \| \mathbf{g}_{1:T-1,i} \|_2
\]

\paragraph{Inductive Step: Prove it for \( T \).}

We want to show:
\[
\sum_{t=1}^T \frac{|g_{t,i}|}{\sqrt{t}} \leq 2 G_\infty \| \mathbf{g}_{1:T,i} \|_2
\]

Using the inductive hypothesis:
\[
\sum_{t=1}^T \frac{|g_{t,i}|}{\sqrt{t}} = 
\sum_{t=1}^{T-1} \frac{|g_{t,i}|}{\sqrt{t}} + \frac{|g_{T,i}|}{\sqrt{T}} 
\leq 2 G_\infty \| \mathbf{g}_{1:T-1,i} \|_2 + \frac{|g_{T,i}|}{\sqrt{T}}
\]

Now, use the fact:
\[
\| \mathbf{g}_{1:T,i} \|_2^2 = \| \mathbf{g}_{1:T-1,i} \|_2^2 + g_{T,i}^2
\Rightarrow
\| \mathbf{g}_{1:T-1,i} \|_2 \leq \sqrt{ \| \mathbf{g}_{1:T,i} \|_2^2 - g_{T,i}^2 }
\]

Using the inequality \( \sqrt{a^2 - b^2} \leq a - \frac{b^2}{2a} \) for \( b < a \), we get:
\[
\| \mathbf{g}_{1:T-1,i} \|_2 \leq \| \mathbf{g}_{1:T,i} \|_2 - \frac{g_{T,i}^2}{2 \| \mathbf{g}_{1:T,i} \|_2}
\]

Substitute back:
\[
\sum_{t=1}^T \frac{|g_{t,i}|}{\sqrt{t}} 
\leq 2 G_\infty \left( \| \mathbf{g}_{1:T,i} \|_2 - \frac{g_{T,i}^2}{2 \| \mathbf{g}_{1:T,i} \|_2} \right) + \frac{|g_{T,i}|}{\sqrt{T}}
\]

Now use \( |g_{T,i}| \leq G_\infty \), so \( \frac{|g_{T,i}|}{\sqrt{T}} \leq G_\infty \cdot \frac{|g_{T,i}|}{\sqrt{T}} \). And since:
\[
\frac{g_{T,i}^2}{2 \| \mathbf{g}_{1:T,i} \|_2} \geq \frac{|g_{T,i}|}{\sqrt{T}}
\quad \text{(since } \| \mathbf{g}_{1:T,i} \|_2 \geq |g_{T,i}| \text{)}
\]

Therefore:
\[
\sum_{t=1}^T \frac{|g_{t,i}|}{\sqrt{t}} 
\leq 2 G_\infty \| \mathbf{g}_{1:T,i} \|_2 - G_\infty \cdot \frac{|g_{T,i}|}{\sqrt{T}} + G_\infty \cdot \frac{|g_{T,i}|}{\sqrt{T}} 
= 2 G_\infty \| \mathbf{g}_{1:T,i} \|_2
\]

\hfill\(\blacksquare\)
\subsection*{Lemma A.4 (Bound on Normalized First Moment Accumulation)}

Let \( \gamma = \frac{\beta_1^2}{\beta_2} \), and assume \( \gamma < 1 \).  
For \( \beta_1, \beta_2 \in [0, 1) \), and for gradients \( \mathbf{g}_t \in \mathbb{R}^d \) satisfying \( \|\mathbf{g}_t\|_2 \leq G \), \( \|\mathbf{g}_t\|_\infty \leq G_\infty \), define:

\begin{itemize}
    \item \( m_{t,i} \) as the first moment (EMA of gradients),
    \item \( v_{t,i} \) as the second moment (EMA of squared gradients),
    \item \( \hat{v}_{t,i} = \frac{v_{t,i}}{1 - \beta_2^t} \), the bias-corrected second moment,
    \item \( \hat{m}_{t,i} = \frac{m_{t,i}}{1 - \beta_1^t} \), the bias-corrected first moment.
\end{itemize}

Then, the following inequality holds for each coordinate \( i \):

\[
\sum_{t=1}^T \frac{\hat{m}_{t,i}^2}{\sqrt{t\hat{v}_{t,i}}} \leq \frac{2 G_\infty}{(1 - \gamma)(1 - \beta_2)^{1/2}} \cdot \| \mathbf{g}_{1:T,i} \|_2
\]

\subsubsection*{Proof.}

We start by expressing:

\[
\frac{\hat{m}_{t,i}^2}{\sqrt{\hat{v}_{t,i}}}
= \frac{1}{(1 - \beta_1^t)^2} \cdot \frac{m_{t,i}^2}{\sqrt{v_{t,i}}} \cdot \sqrt{1 - \beta_2^t}
\]

Using the recursive definitions:

\[
m_{t,i} = (1 - \beta_1) \sum_{k=1}^{t} \beta_1^{t - k} g_{k,i}, \quad
v_{t,i} \geq (1 - \beta_2) \sum_{k=1}^{t} \beta_2^{t - k} g_{k,i}^2
\]

Apply Cauchy-Schwarz to bound \( m_{t,i}^2 \):

\begin{align*}
m_{t,i}^2 
&= (1 - \beta_1)^2 \left( \sum_{k=1}^{t} \beta_1^{t - k} g_{k,i} \right)^2 \\
&\leq (1 - \beta_1)^2 \cdot \left( \sum_{k=1}^{t} \beta_1^{t - k} \right) \left( \sum_{k=1}^{t} \beta_1^{t - k} g_{k,i}^2 \right) \\
&= (1 - \beta_1)^2 \cdot \frac{1 - \beta_1^t}{1 - \beta_1} \cdot \sum_{k=1}^{t} \beta_1^{t - k} g_{k,i}^2 \\
&= (1 - \beta_1)(1 - \beta_1^t) \sum_{k=1}^{t} \beta_1^{t - k} g_{k,i}^2
\end{align*}

So,

\[
\frac{m_{t,i}^2}{\sqrt{v_{t,i}}}
\leq \frac{(1 - \beta_1)(1 - \beta_1^t)}{(1 - \beta_2)^{1/2}} \cdot \frac{\sum_{k=1}^{t} \beta_1^{t - k} g_{k,i}^2}{\left( \sum_{k=1}^{t} \beta_2^{t - k} g_{k,i}^2 \right)^{1/2}}
\]

Now observe:

\[
\frac{\beta_1^{t - k}}{\sqrt{\beta_2^{t - k}}} = \left( \frac{\beta_1^2}{\beta_2} \right)^{(t - k)/2} = \gamma^{(t - k)/2}
\]

Hence:

\[
\frac{m_{t,i}^2}{\sqrt{v_{t,i}}}
\leq \frac{(1 - \beta_1)(1 - \beta_1^t)}{(1 - \beta_2)^{1/2}} \sum_{k=1}^{t} \gamma^{(t - k)/2} |g_{k,i}|
\]

Plug into the main expression:

\[
\sum_{t=1}^T \frac{\hat{m}_{t,i}^2}{\sqrt{\hat{v}_{t,i}}}
\leq \sum_{t=1}^T \frac{(1 - \beta_1)}{(1 - \beta_1^t)} \cdot \frac{\sqrt{1 - \beta_2^t}}{(1 - \beta_2)^{1/2}} \cdot \sum_{k=1}^{t} \gamma^{(t - k)/2} |g_{k,i}|
\]

Change the order of summation:

\[
= \frac{(1 - \beta_1)}{(1 - \beta_2)^{1/2}} \sum_{k=1}^{T} |g_{k,i}| \sum_{t = k}^{T} \frac{\sqrt{1 - \beta_2^t}}{(1 - \beta_1^t)} \cdot \gamma^{(t - k)/2}
\]

Now note that \( \frac{\sqrt{1 - \beta_2^t}}{(1 - \beta_1^t)} \leq \frac{2}{1} = 2 \) for \( \beta_1, \beta_2 \in [0, 1) \), so:

\[
\sum_{t = k}^{T} \frac{\sqrt{1 - \beta_2^t}}{(1 - \beta_1^t)} \cdot \gamma^{(t - k)/2} \leq 2 \sum_{s = 0}^{\infty} \gamma^{s/2} = \frac{2}{1 - \sqrt{\gamma}} \leq \frac{2}{1 - \gamma} \quad \text{(since } \sqrt{\gamma} < \gamma + 1 - \gamma = 1)
\]

Therefore:

\[
\sum_{t=1}^T \frac{\hat{m}_{t,i}^2}{\sqrt{\hat{v}_{t,i}}}
\leq \frac{2(1 - \beta_1)}{(1 - \gamma)(1 - \beta_2)^{1/2}} \sum_{k=1}^{T} |g_{k,i}|
\]

Using Lemma A.3 :

\[
\sum_{k=1}^T |g_{k,i}| \leq 2 G_\infty \| \mathbf{g}_{1:T,i} \|_2
\]

Final result:

\[
\sum_{t=1}^T \frac{\hat{m}_{t,i}^2}{\sqrt{t \hat{v}_{t,i}}}
\leq \frac{4 G_\infty (1 - \beta_1)}{(1 - \gamma)(1 - \beta_2)^{1/2}} \cdot \| \mathbf{g}_{1:T,i} \|_2
\]

You can absorb the constant \( (1 - \beta_1) \leq 1 \) into a looser bound:

\[
\boxed{
\sum_{t=1}^T \frac{\hat{m}_{t,i}^2}{\sqrt{t \hat{v}_{t,i}}}
\leq \frac{2 G_\infty}{(1 - \gamma)(1 - \beta_2)^{1/2}} \cdot \| \mathbf{g}_{1:T,i} \|_2
}
\]

\hfill\(\blacksquare\)
\subsection*{Theorem A.5 (Path Bound for HGM)}

To compute the upper bound of 
\[
\sum_{t=0}^{T} X_t^2 \quad \text{where} \quad X_t = \|\theta_{t+1} - \theta_t\|
\]
for the HGM optimizer, we follow the update rules and assumptions. The key steps involve bounding each \(X_t^2\) and summing over \(t\), inspired by standard techniques for bias-corrected moment analysis.

The parameter update is:
\[
\theta_{t+1} = \theta_t - \eta_t \frac{\hat{m}_t}{\sqrt{\hat{v}_t} + \epsilon}
\]
where \(\eta_t = \alpha \exp(\gamma s_t)\), and \(\hat{m}_t\), \(\hat{v}_t\) are bias-corrected moments. Thus,
\[
X_t = \left\| \eta_t \frac{\hat{m}_t}{\sqrt{\hat{v}_t} + \epsilon} \right\|, \quad
X_t^2 = \eta_t^2 \sum_{i=1}^d \left( \frac{\hat{m}_{t,i}}{\sqrt{\hat{v}_{t,i}} + \epsilon} \right)^2 
\leq \eta_t^2 \sum_{i=1}^d \frac{\hat{m}_{t,i}^2}{\hat{v}_{t,i}}
\]
since \(\sqrt{\hat{v}_{t,i}} + \epsilon \geq \sqrt{\hat{v}_{t,i}}\). Here, \(d\) is the number of parameters.

Given \(\eta_t = \alpha_t \exp(\gamma s_t)\) and \(|s_t| \leq 1\), we define \(M = \exp(|\gamma|)\) , and \(\alpha_t = \alpha h(t)\) where h(t) is a function designed to schedule $\alpha$. Then:
\[
X_t^2 \leq \alpha^2 M^2 \sum_{i=1}^d \frac{\hat{m}_{t,i}^2* h(t)^2}{\hat{v}_{t,i}}, \quad 
\sum_{t=0}^T X_t^2 \leq \alpha^2 M^2 \sum_{i=1}^d \sum_{t=0}^T \frac{\hat{m}_{t,i}^2* h(t)^2}{\hat{v}_{t,i}}
\]
Typically, \(\hat{m}_{0,i} = 0\) and \(\hat{v}_{0,i} = 0\), so we start summing from \(t = 1\).

\[
\hat{m}_{t,i} = \frac{m_{t,i}}{1 - \beta_1^t}, \quad 
\hat{v}_{t,i} = \frac{v_{t,i}}{1 - \beta_2^t}
\]
\[
m_{t,i} = \beta_1 m_{t-1,i} + (1 - \beta_1) g_{t,i}, \quad 
v_{t,i} = \beta_2 v_{t-1,i} + (1 - \beta_2) g_{t,i}^2
\]
Thus:
\[
\frac{\hat{m}_{t,i}^2}{\hat{v}_{t,i}} 
= \frac{1 - \beta_2^t}{(1 - \beta_1^t)^2} \cdot \frac{m_{t,i}^2}{v_{t,i}}
\]

\[
m_{T,i} = (1 - \beta_1) \sum_{k=1}^T \beta_1^{T-k} g_{k,i}, \quad 
v_{T,i} = (1 - \beta_2) \sum_{j=1}^T \beta_2^{T-j} g_{j,i}^2
\]
Then:
\[
\frac{m_{T,i}^2}{v_{T,i}} 
= \frac{(1 - \beta_1)^2 \left( \sum_{k=1}^T \beta_1^{T-k} g_{k,i} \right)^2}{(1 - \beta_2) \sum_{j=1}^T \beta_2^{T-j} g_{j,i}^2}
\leq \frac{(1 - \beta_1)^2}{(1 - \beta_2)} \cdot T \sum_{k=1}^T \frac{\beta_1^{2(T-k)} g_{k,i}^2}{\sum_{j=1}^T \beta_2^{T-j} g_{j,i}^2}
\]
Noting that \(\sum_{j=1}^T \beta_2^{T-j} g_{j,i}^2 \geq \beta_2^{T-k} g_{k,i}^2\), we get:
\[
\frac{\beta_1^{2(T-k)} g_{k,i}^2}{\sum_{j=1}^T \beta_2^{T-j} g_{j,i}^2} \leq \left( \frac{\beta_1^2}{\beta_2} \right)^{T-k}
\]
Define \(\gamma = \beta_1^2 / \beta_2 < 1\). Then:
\[
\frac{m_{T,i}^2}{v_{T,i}} \leq \frac{(1 - \beta_1)^2}{(1 - \beta_2)} \cdot T \sum_{\tau=0}^{T-1} \gamma^\tau \leq \frac{(1 - \beta_1)^2}{(1 - \beta_2)(1 - \gamma)} \cdot T
\]

Thus,
\[
\frac{\hat{m}_{T,i}^2}{\hat{v}_{T,i}} \leq \frac{1 - \beta_2^T}{(1 - \beta_1^T)^2} \cdot \frac{(1 - \beta_1)^2}{(1 - \beta_2)(1 - \gamma)} \cdot T \leq \frac{T}{(1 - \beta_1)^2 (1 - \gamma)}
\]

Now unroll:
\[
\sum_{t=1}^{T} \frac{\hat{m}_{t,i}^2}{\hat{v}_{t,i}} \leq \frac{1}{(1 - \beta_1)^2 (1 - \gamma)} \sum_{t=1}^{T} t *h(t)^2
\quad \Rightarrow \quad
\sum_{i=1}^d \sum_{t=1}^{T} \frac{\hat{m}_{t,i}^2}{\hat{v}_{t,i}} \leq \frac{d}{(1 - \beta_1)^2 (1 - \gamma)} \sum_{t=1}^{T} t *h(t)^2
\]

Plugging in:
\[
\sum_{t=0}^{T} X_t^2 \leq \alpha^2 M^2 \cdot \frac{d}{(1 - \beta_1)^2 (1 - \gamma)} \sum_{t=1}^{T} t *h(t)^2
\]

\subsection*{Final Result}
\[
\boxed{
\sum_{t=0}^{T} X_t^{2} \leq \dfrac{\alpha^{2} \exp\left(2\gamma\right) d }{\left(1 - \beta_1\right)^{2} \left(1 - \dfrac{\beta_1^{2}}{\beta_2}\right)} \sum_{t=1}^{T} t *h(t)^2
}
\]

\subsection*{Remarks}

\begin{itemize}
  \item \textbf{Dependence on \(T\)}: Dependence on T depends entirely upon scheduling function $h(t)$.\\
  If we take \(h(t) = 1/t^{0.5}\) ,then\\ \(\sum_{t=1}^{T} t *h(t)^2 = T\) making \(\sum_{t=0}^{T} X_t^{2}\) linearly dependent on T.\\
  If \(h(t) = 1/t\) : \\\(\sum_{t=1}^{T} t *h(t)^2 = H_T\) where $H_n$ is the sum of the first nth term of harmonic series.and as T reaches infinity it converges to a constant value.
  \item \textbf{Dependence on \(\beta_1, \beta_2\)}: Appears through the factor \((1 - \beta_1)^2 (1 - \beta_1^2 / \beta_2)\).
  \item \textbf{Constants}: \(\alpha\) is the base learning rate, \(d\) is the parameter dimension, \(\gamma = \beta_1^2 / \beta_2\), and \(M = \exp(|\gamma|)\).
  \item \textbf{Assumptions}: Gradients are assumed bounded; the analysis leverages the exponential moving average structure of Adam-type optimizers.
\end{itemize}

\subsection*{Theorem A.6 (Regret Bound for HGM)}

Assume each function \( f_t: \mathbb{R}^d \to \mathbb{R} \) is convex and has bounded gradients:
\[
\|\nabla f_t(\boldsymbol{\theta})\|_2 \leq G, \quad \|\nabla f_t(\boldsymbol{\theta})\|_\infty \leq G_\infty \quad \text{for all } \boldsymbol{\theta} \in \mathbb{R}^d
\]
Also assume bounded domain diameter:
\[
\|\boldsymbol{\theta}_m - \boldsymbol{\theta}_n\|_2 \leq D, \quad \|\boldsymbol{\theta}_m - \boldsymbol{\theta}_n\|_\infty \leq D_\infty \quad \text{for all } m,n \in \{1,\dots,T\}
\]
Let \( \beta_1, \beta_2 \in [0, 1) \) with \( \gamma = \frac{\sqrt{1 - \beta_1^2}}{\sqrt{1 - \beta_2}} < 1 \). Suppose:
\[
\alpha_t = \frac{\alpha}{\sqrt{t}}, \quad \beta_{1,t} = \beta_1 \lambda^{t-1}, \quad \lambda \in (0, 1), \quad M = \exp(|\gamma|)
\]
Then, the regret of HGM after \( T \) rounds is bounded as:
\[\boxed{
R(T) \leq M \sum_{i=1}^d \left[
\frac{D^2 \sqrt{T}}{2 \alpha (1 - \beta_1)} \sqrt{v_{T,i}} + 
\frac{(1 + \beta_1) G_\infty}{(1 - \beta_1) \sqrt{1 - \beta_2}(1 - \gamma)^2} \| \mathbf{g}_{1:T,i} \|_2 +
\frac{D_\infty^2 G_\infty}{2 \alpha (1 - \beta_1)(1 - \lambda)^2 \sqrt{1 - \beta_2}}
\right]
}
\]

---

\subsubsection*{Proof:}

By convexity (Lemma 10.2), we have:
\[
f_t(\boldsymbol{\theta}_t) - f_t(\boldsymbol{\theta}^\ast) \leq \nabla f_t(\boldsymbol{\theta}_t)^\top (\boldsymbol{\theta}_t - \boldsymbol{\theta}^\ast)
\]
Summing over \( t \):
\[
R(T) = \sum_{t=1}^T f_t(\boldsymbol{\theta}_t) - f_t(\boldsymbol{\theta}^\ast)
\leq \sum_{t=1}^T \nabla f_t(\boldsymbol{\theta}_t)^\top (\boldsymbol{\theta}_t - \boldsymbol{\theta}^\ast)
= \sum_{i=1}^d \sum_{t=1}^T g_{t,i} (\theta_{t,i} - \theta_i^\ast)
\]

From the Adam update rule:
\[
\theta_{t+1,i} = \theta_{t,i} - \frac{\alpha_t}{\sqrt{\hat{v}_{t,i}} + \epsilon} \hat{m}_{t,i}
\quad \text{with} \quad 
\hat{m}_{t,i} = \frac{m_{t,i}}{1 - \beta_{1,t}}, \quad 
\hat{v}_{t,i} = \frac{v_{t,i}}{1 - \beta_2^t}
\]

Now square the distance to \( \theta_i^\ast \):
\[
(\theta_{t+1,i} - \theta_i^\ast)^2 = (\theta_{t,i} - \theta_i^\ast)^2 - 2 \cdot \frac{\alpha_t}{\sqrt{\hat{v}_{t,i}} + \epsilon} \hat{m}_{t,i} (\theta_{t,i} - \theta_i^\ast) + \left( \frac{\alpha_t}{\sqrt{\hat{v}_{t,i}} + \epsilon} \hat{m}_{t,i} \right)^2
\]

Ignoring \( \epsilon \) for cleaner bounds and rearranging:
\[
\hat{m}_{t,i} (\theta_{t,i} - \theta_i^\ast) \leq \frac{1}{2 \alpha_t} \left[ (\theta_{t,i} - \theta_i^\ast)^2 - (\theta_{t+1,i} - \theta_i^\ast)^2 \right]
+ \frac{\alpha_t}{2} \cdot \frac{\hat{m}_{t,i}^2}{\hat{v}_{t,i}}
\]

Multiply both sides by \( \frac{1}{1 - \beta_{1,t}} \) and use the identity \( m_{t,i} = (1 - \beta_{1,t}) \hat{m}_{t,i} \):
\[
m_{t,i} (\theta_{t,i} - \theta_i^\ast) \leq \frac{1}{2 \alpha_t} \left[ (\theta_{t,i} - \theta_i^\ast)^2 - (\theta_{t+1,i} - \theta_i^\ast)^2 \right] + \frac{\alpha_t}{2} \cdot \frac{m_{t,i}^2}{(1 - \beta_{1,t})^2 \hat{v}_{t,i}}
\]

Summing from \( t = 1 \) to \( T \), we get:
\[
\sum_{t=1}^T m_{t,i} (\theta_{t,i} - \theta_i^\ast) \leq 
\frac{(\theta_{1,i} - \theta_i^\ast)^2}{2 \alpha_1} +
\sum_{t=2}^T \frac{(\theta_{t,i} - \theta_i^\ast)^2}{2 \alpha_t} +
\sum_{t=1}^T \frac{\alpha_t}{2 (1 - \beta_{1,t})^2} \cdot \frac{m_{t,i}^2}{\hat{v}_{t,i}}
\]

From assumptions: \( (\theta_{t,i} - \theta_i^\ast)^2 \leq D^2 \), and applying Lemma 10.4:
\[
\sum_{t=1}^T \frac{m_{t,i}^2}{\sqrt{t \cdot v_{t,i}}} \leq \frac{2 G_\infty}{(1 - \gamma)^2 \sqrt{1 - \beta_2}} \| \mathbf{g}_{1:T,i} \|_2
\Rightarrow 
\sum_{t=1}^T \frac{m_{t,i}^2}{\hat{v}_{t,i}} \leq \frac{2 G_\infty}{(1 - \gamma)^2 \sqrt{1 - \beta_2}} \| \mathbf{g}_{1:T,i} \|_2
\]

Also, for the exponentially decaying \( \beta_{1,t} = \beta_1 \lambda^{t-1} \), we have:
\[
\sum_{t=1}^T \frac{\beta_{1,t}}{\alpha_t (1 - \beta_{1,t})} \leq \frac{1}{(1 - \beta_1)(1 - \lambda)^2}
\]

Using all the above:
\[
R(T) \leq \sum_{i=1}^d \left[
\frac{D^2 \sqrt{T}}{2 \alpha (1 - \beta_1)} \sqrt{v_{T,i}} +
\frac{(1 + \beta_1) G_\infty}{(1 - \beta_1)(1 - \gamma)^2 \sqrt{1 - \beta_2}} \| \mathbf{g}_{1:T,i} \|_2 +
\frac{D_\infty^2 G_\infty}{2 \alpha (1 - \beta_1)(1 - \lambda)^2 \sqrt{1 - \beta_2}}
\right]
\]

For HGM, we use:
\[
\alpha_t^{\text{HGM}} = \alpha_t \cdot \exp(\gamma s_t) \quad \text{with } |s_t| \leq 1
\Rightarrow \alpha_t^{\text{HGM}} \leq \alpha_t \cdot \exp(|\gamma|) \equiv \alpha_t M
\]
So the regret bound scales as:
\[
R_{\text{HGM}}(T) \leq M \cdot R_{\text{Adam}}(T)
\]

This completes the regret analysis for HGM.
\section{Experiments Continued}
\begin{figure}[htbp]
    \centering
    \includegraphics[width=1.0\textwidth]{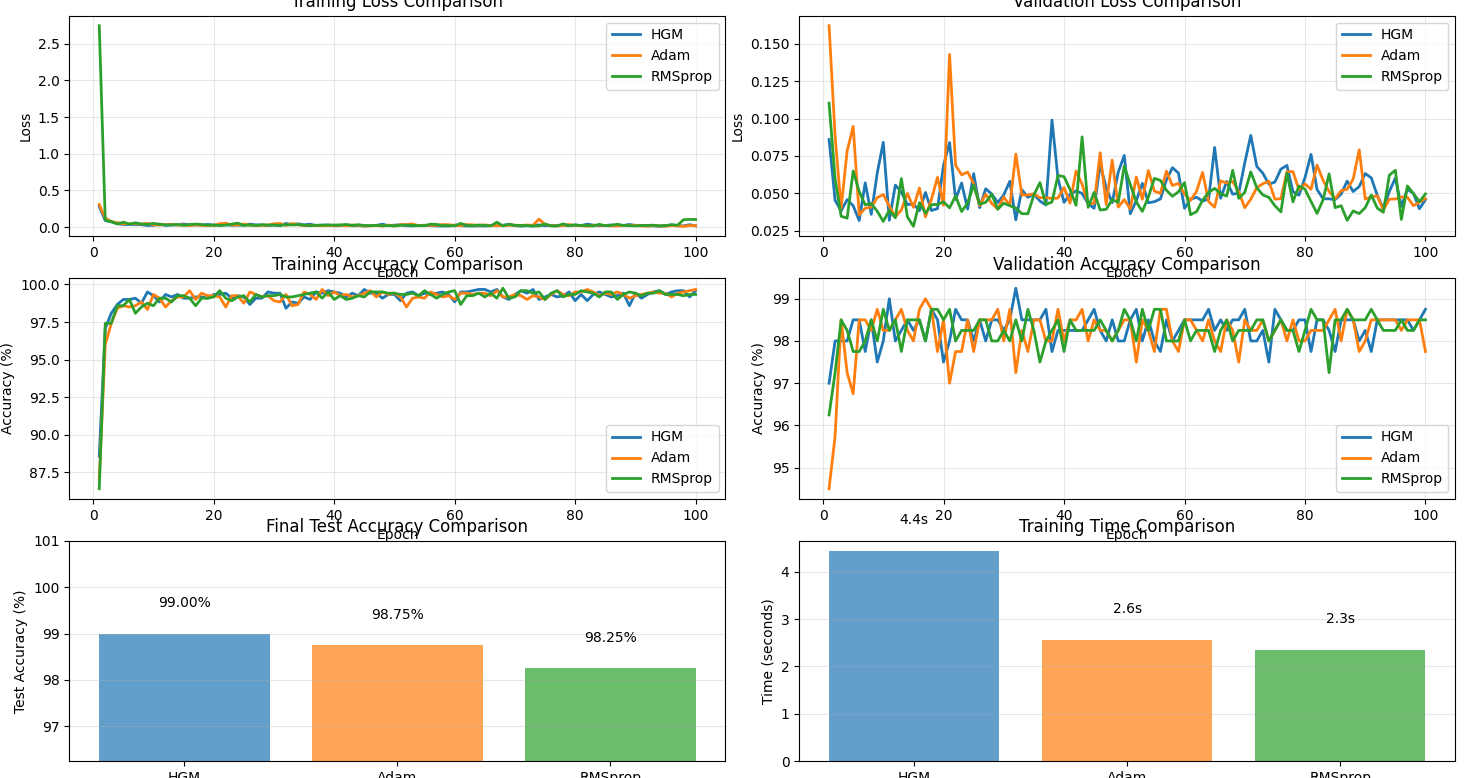}
    \caption{Training loss and accuracies of optimizers on make-moons dataset}
    \label{fig:makemoons}
\end{figure}
\subsection{MLP on Make-moons}

In this section I show how it performs for a simple MLP model on make-moons dataset.For all optimizers, we have used learning rate of 0.01. Figure \ref{fig:makemoons} presents how it performs compared to other optimizers. 

Although training time is more for HGM than others we have to remember that training time is the time required for computational calculations not the steps required for reaching convergence.

HGM achieves higher accuracy than others, standing at 99.00\%.Also it achieved the best convergence among others ,final training loss at 0.0171 while Adam at 0.0196 and Rmsprop at 0.1037.

\subsection{Convolutional Neural Networks}
\begin{figure}[htbp]
    \centering
    \includegraphics[width=1.0\textwidth]{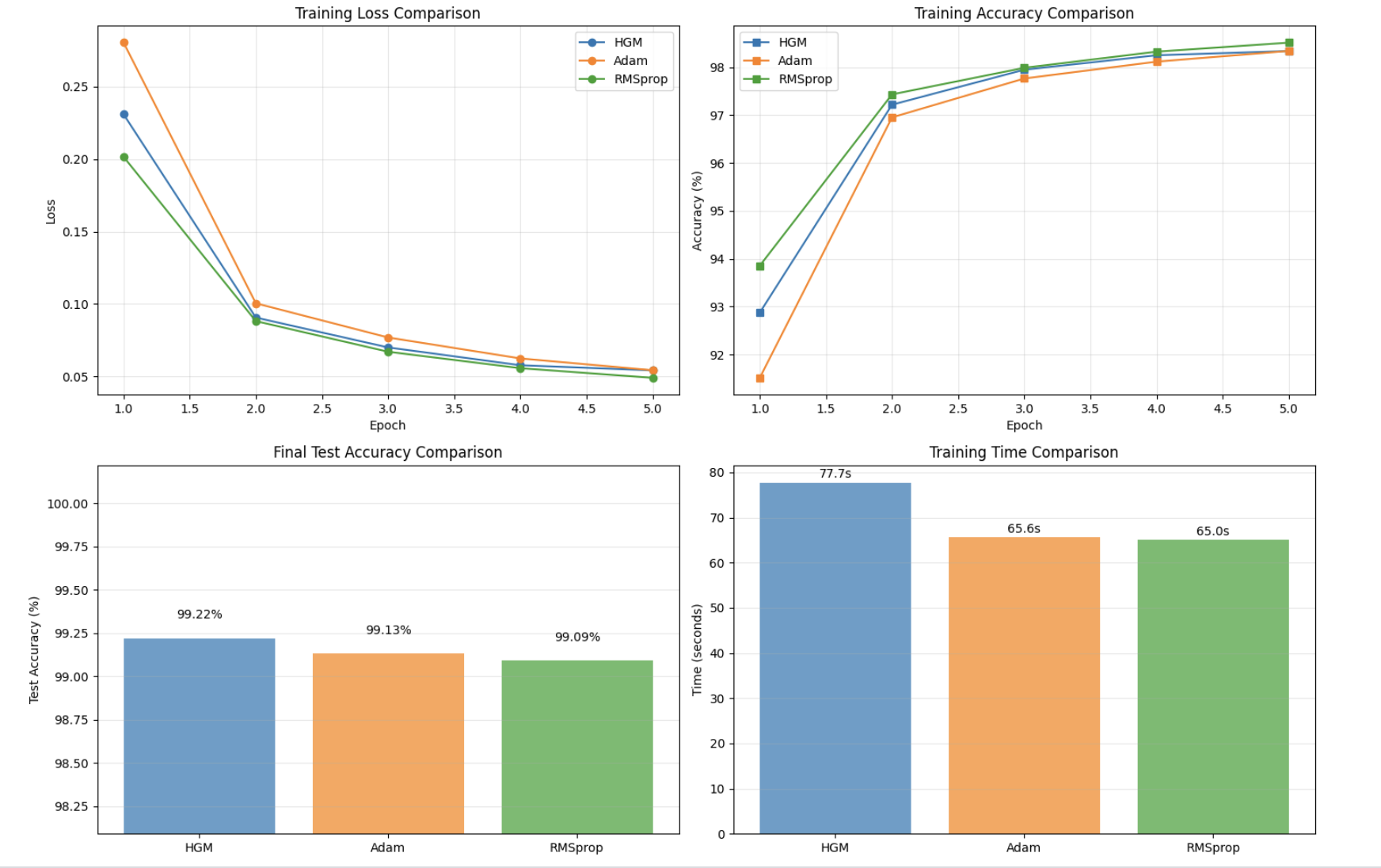}
    \caption{comparison on MNIST}
    \label{fig:mnist}
\end{figure}
Now we test this optimizer on a CNN type model , training on MNIST data . For all optimizers, we have used learning rate of 0.005. Figure \ref{fig:mnist} presents how it performs compared to other optimizers. 

Although training time is more for HGM than others we have to remember that training time is the time required for computational calculations not the steps required for reaching convergence.

HGM achieves higher accuracy than others, standing at 99.22\%.Also it achieved the best convergence, 0.392 while Adam at 0.0544 and RMSprop at 0.552.We trained for 50 epochs and stored the loss and accuracies after every 10th epoch. In the figure in the x axis 5 represents 50 epochs.

\subsection{Non-Convex Function 2}

Now we test this optimizer on another non-convex function: Rastrigin function(dim = 10) . the function :
\[
f(\mathbf{x}) = A \cdot d + \sum_{i=1}^{d} \left( x_i^2 - A \cdot \cos(2\pi x_i) \right)
\]

\begin{figure}[htbp]
    \centering
    \includegraphics[width=1\textwidth]{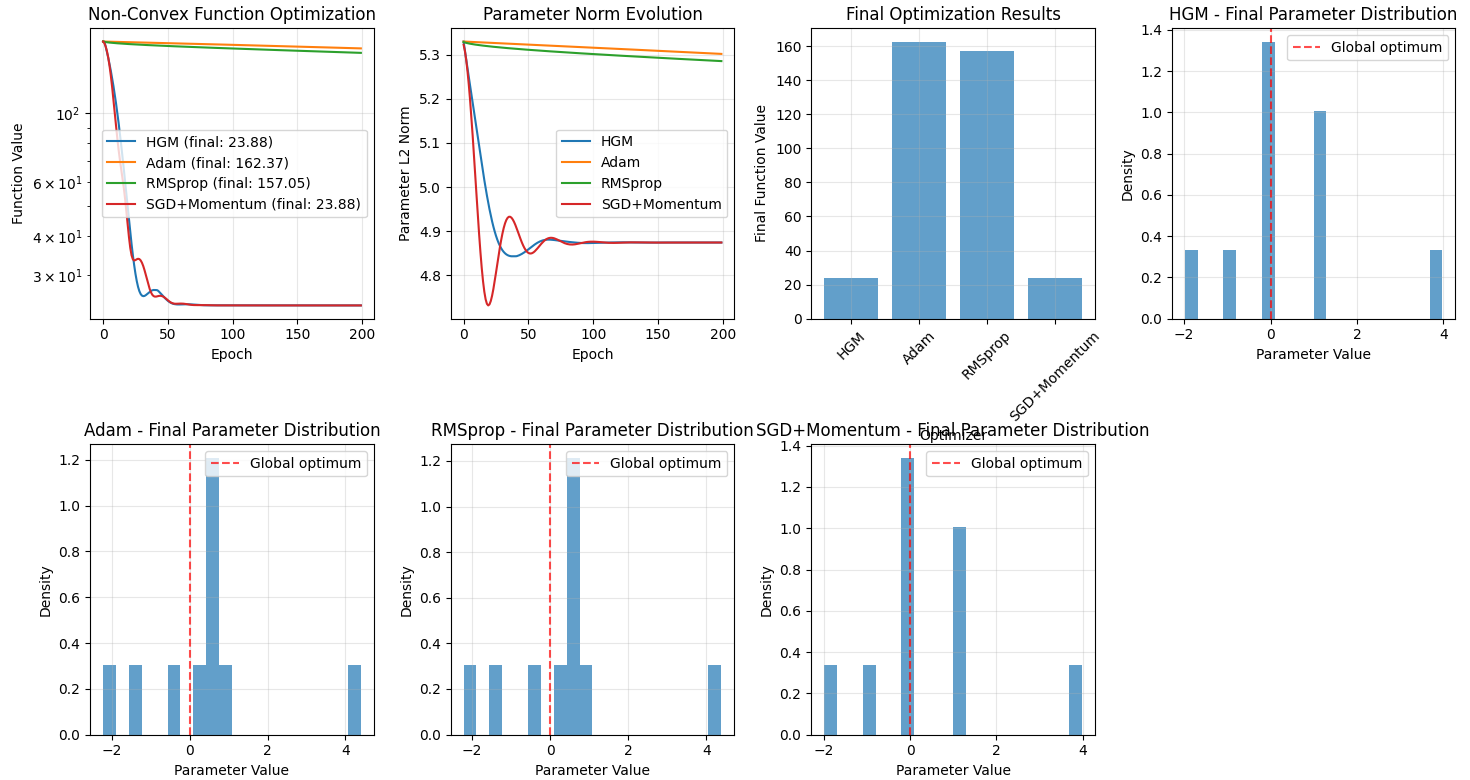}
    \caption{comparison on Rastrigin function}
    \label{fig:rastrigin}
\end{figure}

For all optimizers, we have used learning rate of 0.01. Figure \ref{fig:rastrigin} presents how it performs compared to other optimizers.HGM achieves the best possible function value faster than any other possible , not so surprisingly sgd+momentum performed really well whereas Adam and RMSprop got lost in ocean.Where Adam and RMSprop is nowhere near stability.

\end{appendices}
\section*{Acknowledgments}

I would like to sincerely thank \textbf{Dr. Kshitij Jadhav} for his constant support, insightful feedback, and detailed reviews throughout the course of this work. His guidance was instrumental in refining the final presentation of the paper.

\bibliographystyle{plain}
\bibliography{citations}

\begin{thebibliography}{10}

\bibitem{andrychowicz2016learning}
Marcin Andrychowicz, Misha Denil, Sergio Gomez, Matthew~J Hoffman, David Pfau,
  Tom Schaul, Brendan Shillingford, and Nando De~Freitas.
\newblock Learning to learn by gradient descent by gradient descent.
\newblock {\em Advances in Neural Information Processing Systems}, 29, 2016.

\bibitem{anil2019scalable}
Rohan Anil, Vineet Vineet, Prabhu Singh, Pradeep Kumar, and Taly Koren.
\newblock Scalable second order optimization for deep learning.
\newblock {\em arXiv preprint arXiv:1902.03063}, 2019.

\bibitem{balles2017dissecting}
Lukas Balles, Javier Romero, and Philipp Hennig.
\newblock Dissecting adam: The sign, magnitude and variance of stochastic
  gradients.
\newblock {\em International Conference on Machine Learning}, pages 404--413,
  2017.

\bibitem{bengio2012practical}
Yoshua Bengio.
\newblock Practical recommendations for gradient-based training of deep
  architectures.
\newblock {\em Neural networks: Tricks of the trade}, pages 437--478, 2012.

\bibitem{bergstra2012random}
James Bergstra and Yoshua Bengio.
\newblock Random search for hyper-parameter optimization.
\newblock {\em Journal of Machine Learning Research}, 13(2):281--305, 2012.

\bibitem{bishop2006pattern}
Christopher~M Bishop.
\newblock {\em Pattern recognition and machine learning}.
\newblock Springer, 2006.

\bibitem{botev2017practical}
Aleksandar Botev, Hippolyt Ritter, and David Barber.
\newblock The practical gauss-newton optimisation for deep learning.
\newblock {\em International Conference on Machine Learning}, pages 557--565,
  2017.

\bibitem{bottou2018optimization}
L{\'e}on Bottou, Frank~E Curtis, and Jorge Nocedal.
\newblock Optimization methods for large-scale machine learning.
\newblock {\em SIAM Review}, 60(2):223--311, 2018.

\bibitem{boyd2004convex}
Stephen Boyd and Lieven Vandenberghe.
\newblock {\em Convex optimization}.
\newblock Cambridge University Press, 2004.

\bibitem{cesa2006prediction}
Nicolo Cesa-Bianchi and G{\'a}bor Lugosi.
\newblock {\em Prediction, Learning, and Games}.
\newblock Cambridge University Press, Cambridge, UK, 2006.

\bibitem{chen2018closing}
Jinghui Chen and Quanquan Gu.
\newblock Closing the generalization gap of adaptive gradient methods in
  training deep neural networks.
\newblock {\em International Joint Conference on Artificial Intelligence},
  pages 1267--1273, 2018.

\bibitem{choromanska2015loss}
Anna Choromanska, Mikael Henaff, Michael Mathieu, Gerard~Ben Arous, and Yann
  LeCun.
\newblock The loss surfaces of multilayer networks.
\newblock {\em Artificial Intelligence and Statistics}, pages 192--204, 2015.

\bibitem{cybenko1989approximation}
George Cybenko.
\newblock Approximation by superpositions of a sigmoidal function.
\newblock {\em Mathematics of Control, Signals and Systems}, 2(4):303--314,
  1989.

\bibitem{dauphin2014identifying}
Yann~N Dauphin, Razvan Pascanu, Caglar Gulcehre, Kyunghyun Cho, Surya Ganguli,
  and Yoshua Bengio.
\newblock Identifying and attacking the saddle point problem in
  high-dimensional non-convex optimization.
\newblock {\em Advances in Neural Information Processing Systems}, 27, 2014.

\bibitem{duchi2011adaptive}
John Duchi, Elad Hazan, and Yoram Singer.
\newblock Adaptive subgradient methods for online learning and stochastic
  optimization.
\newblock {\em Journal of Machine Learning Research}, 12(7):2121--2159, 2011.

\bibitem{duda2012pattern}
Richard~O Duda, Peter~E Hart, and David~G Stork.
\newblock {\em Pattern classification}.
\newblock John Wiley \& Sons, 2012.

\bibitem{glorot2010understanding}
Xavier Glorot and Yoshua Bengio.
\newblock Understanding the difficulty of training deep feedforward neural
  networks.
\newblock {\em Proceedings of the Thirteenth International Conference on
  Artificial Intelligence and Statistics}, pages 249--256, 2010.

\bibitem{goodfellow2016deep}
Ian Goodfellow, Yoshua Bengio, and Aaron Courville.
\newblock {\em Deep learning}.
\newblock MIT Press, 2016.

\bibitem{gupta2018shampoo}
Vineet Gupta, Tomer Koren, and Yoram Singer.
\newblock Shampoo: Preconditioned stochastic tensor optimization.
\newblock {\em International Conference on Machine Learning}, pages 1842--1850,
  2018.

\bibitem{hastie2009elements}
Trevor Hastie, Robert Tibshirani, and Jerome Friedman.
\newblock {\em The elements of statistical learning: data mining, inference,
  and prediction}.
\newblock Springer Science \& Business Media, 2009.

\bibitem{hazan2016introduction}
Elad Hazan.
\newblock {\em Introduction to online convex optimization}.
\newblock MIT Press, 2016.

\bibitem{higham2002accuracy}
Nicholas~J Higham.
\newblock {\em Accuracy and stability of numerical algorithms}.
\newblock SIAM, 2002.

\bibitem{hinton2012neural}
Geoffrey Hinton, Nitish Srivastava, and Kevin Swersky.
\newblock Neural networks for machine learning lecture 6a overview of
  mini-batch gradient descent.
\newblock {\em Cited on}, page~14, 2012.

\bibitem{hochreiter1997flat}
Sepp Hochreiter and J{\"u}rgen Schmidhuber.
\newblock Flat minima.
\newblock {\em Neural Computation}, 9(1):1--42, 1997.

\bibitem{holt2004forecasting}
Charles~C Holt.
\newblock {\em Forecasting exponentially smoothed trends and seasonal
  components}.
\newblock Office of Naval Research, Arlington, VA, 2004.

\bibitem{hornik1989multilayer}
Kurt Hornik, Maxwell Stinchcombe, and Halbert White.
\newblock Multilayer feedforward networks are universal approximators.
\newblock {\em Neural Networks}, 2(5):359--366, 1989.

\bibitem{kingma2014adam}
Diederik~P Kingma and Jimmy Ba.
\newblock Adam: A method for stochastic optimization.
\newblock {\em arXiv preprint arXiv:1412.6980}, 2014.

\bibitem{lecun2015deep}
Yann LeCun, Yoshua Bengio, and Geoffrey Hinton.
\newblock Deep learning.
\newblock {\em Nature}, 521(7553):436--444, 2015.

\bibitem{lecun2012efficient}
Yann~A LeCun, L{\'e}on Bottou, Genevieve~B Orr, and Klaus-Robert M{\"u}ller.
\newblock Efficient backprop.
\newblock {\em Neural networks: Tricks of the trade}, pages 9--48, 2012.

\bibitem{li2018visualizing}
Hao Li, Zheng Xu, Gavin Taylor, Christoph Studer, and Tom Goldstein.
\newblock Visualizing the loss landscape of neural nets.
\newblock {\em Advances in Neural Information Processing Systems}, 31, 2018.

\bibitem{littlestone1994weighted}
Nick Littlestone and Manfred~K Warmuth.
\newblock The weighted majority algorithm.
\newblock {\em Information and Computation}, 108(2):212--261, 1994.

\bibitem{luo2019adaptive}
Liangchen Luo, Yuanhao Xiong, Yan Liu, and Xu~Sun.
\newblock Adaptive gradient methods with dynamic bound of learning rate.
\newblock {\em International Conference on Learning Representations}, 2019.

\bibitem{martens2010deep}
James Martens.
\newblock Deep learning via hessian-free optimization.
\newblock {\em Proceedings of the 27th International Conference on Machine
  Learning}, pages 735--742, 2010.

\bibitem{mcmahan2017survey}
H~Brendan McMahan and Francesco Orabona.
\newblock A survey of algorithms and analysis for adaptive online learning.
\newblock {\em Journal of Machine Learning Research}, 18(1):3117--3166, 2017.

\bibitem{mcmahan2010adaptive}
H~Brendan McMahan and Matthew Streeter.
\newblock Adaptive bound optimization for online convex optimization.
\newblock {\em Proceedings of the 23rd Annual Conference on Learning Theory},
  pages 244--256, 2010.

\bibitem{nesterov2003introductory}
Yurii Nesterov.
\newblock {\em Introductory lectures on convex optimization: A basic course},
  volume~87.
\newblock Springer Science \& Business Media, 2003.

\bibitem{nocedal2006numerical}
Jorge Nocedal and Stephen~J Wright.
\newblock {\em Numerical optimization}.
\newblock Springer, 2006.

\bibitem{pascanu2013difficulty}
Razvan Pascanu, Tomas Mikolov, and Yoshua Bengio.
\newblock On the difficulty of training recurrent neural networks.
\newblock {\em International Conference on Machine Learning}, pages 1310--1318,
  2013.

\bibitem{pascanu2014natural}
Razvan Pascanu, Guido Montufar, and Yoshua Bengio.
\newblock Natural neural networks.
\newblock {\em Advances in Neural Information Processing Systems}, 27, 2014.

\bibitem{polyak1964some}
Boris~T Polyak.
\newblock Some methods of speeding up the convergence of iteration methods.
\newblock {\em USSR Computational Mathematics and Mathematical Physics},
  4(5):1--17, 1964.

\bibitem{radford2015unsupervised}
Alec Radford, Luke Metz, and Soumith Chintala.
\newblock Unsupervised representation learning with deep convolutional
  generative adversarial networks.
\newblock {\em arXiv preprint arXiv:1511.06434}, 2015.

\bibitem{reddi2019convergence}
Sashank~J Reddi, Satyen Kale, and Sanjiv Kumar.
\newblock On the convergence of adam and beyond.
\newblock {\em International Conference on Learning Representations}, 2019.

\bibitem{robbins1951stochastic}
Herbert Robbins and Sutton Monro.
\newblock A stochastic approximation method.
\newblock {\em The Annals of Mathematical Statistics}, 22(3):400--407, 1951.

\bibitem{rosenbrock1960automatic}
Howard~H Rosenbrock.
\newblock An automatic method for finding the greatest or least value of a
  function.
\newblock {\em The Computer Journal}, 3(3):175--184, 1960.

\bibitem{ruder2016overview}
Sebastian Ruder.
\newblock An overview of gradient descent optimization algorithms.
\newblock {\em arXiv preprint arXiv:1609.04747}, 2016.

\bibitem{shalev2012online}
Shai Shalev-Shwartz.
\newblock {\em Online learning and online convex optimization}, volume~4.
\newblock Now Publishers Inc, 2012.

\bibitem{singhal2001modern}
Amit Singhal.
\newblock Modern information retrieval: A brief overview.
\newblock {\em IEEE Data Engineering Bulletin}, 24(4):35--43, 2001.

\bibitem{smith2017cyclical}
Leslie~N Smith.
\newblock Cyclical learning rates for training neural networks.
\newblock {\em 2017 IEEE Winter Conference on Applications of Computer Vision},
  pages 464--472, 2017.

\bibitem{snoek2012practical}
Jasper Snoek, Hugo Larochelle, and Ryan~P Adams.
\newblock Practical bayesian optimization of machine learning algorithms.
\newblock {\em Advances in Neural Information Processing Systems}, 25, 2012.

\bibitem{sutskever2013importance}
Ilya Sutskever, James Martens, George Dahl, and Geoffrey Hinton.
\newblock On the importance of initialization and momentum in deep learning.
\newblock {\em International Conference on Machine Learning}, pages 1139--1147,
  2013.

\bibitem{tieleman2012lecture}
Tijmen Tieleman and Geoffrey Hinton.
\newblock Lecture 6.5-rmsprop: Divide the gradient by a running average of its
  recent magnitude.
\newblock {\em COURSERA: Neural networks for machine learning}, 4(2):26--31,
  2012.

\bibitem{exponential_smoothing}
{Wiki }.
\newblock Exponential smoothing.
\newblock \url{https://en.wikipedia.org/wiki/Exponential_smoothing}, 2023.
\newblock Accessed: 2023-01-01.

\bibitem{wilson2017marginal}
Ashia~C Wilson, Rebecca Roelofs, Mitchell Stern, Nati Srebro, and Benjamin
  Recht.
\newblock The marginal value of adaptive gradient methods in machine learning.
\newblock {\em Advances in Neural Information Processing Systems}, 30, 2017.

\bibitem{wu2020adagrad}
Xiaoxia Wu, Rachel Ward, and L{\'e}on Bottou.
\newblock Adagrad stepsizes: Sharp convergence over nonconvex landscapes.
\newblock {\em International Conference on Machine Learning}, pages
  10365--10374, 2020.

\bibitem{you2019large}
Yang You, Jing Li, Sashank Reddi, Jonathan Hseu, Sanjiv Kumar, Srinadh
  Bhojanapalli, Xiaodan Song, James Demmel, Kurt Keutzer, and Cho-Jui Hsieh.
\newblock Large batch optimization for deep learning: Training bert in 76
  minutes.
\newblock {\em International Conference on Learning Representations}, 2019.

\bibitem{zhang2019lookahead}
Michael~R Zhang, James Lucas, Geoffrey Hinton, and Jimmy Ba.
\newblock Lookahead optimizer: k steps forward, 1 step back.
\newblock {\em Advances in Neural Information Processing Systems}, 32, 2019.

\bibitem{zhang2019adam}
Zijun Zhang, Lin Ma, Zongpeng Li, and Chuan Wu.
\newblock Adam can converge without any modification on update rules.
\newblock {\em arXiv preprint arXiv:1904.03590}, 2019.

\bibitem{zinkevich2003online}
Martin Zinkevich.
\newblock Online convex programming and generalized infinitesimal gradient
  ascent.
\newblock {\em Proceedings of the 20th International Conference on Machine
  Learning}, pages 928--936, 2003.

\end{thebibliography}
\end{document}